\documentclass[11pt]{article}

\usepackage{amscd,amsmath, amssymb, fancyhdr, graphicx, url}

\numberwithin{equation}{section}

\newcommand{\version}{version 1.0,\ \ Apr 08, 2016}

\makeatletter
\def\x@arrow{\DOTSB\Relbar}
\def\xlongrightarrowfill@{\arrowfill@\relbar\relbar\longrightarrow}
\newcommand{\xlongrightarrow}[2][]{%
        \ext@arrow 0099\xlongrightarrowfill@{#1}{#2}}
\makeatother

\def\eqref#1{(\ref{#1})}

\newcommand{\Z}{{\Bbb Z}}
\def\C{{\Bbb C}}

\newcommand{\R}{{\Bbb R}}
\newcommand{\Q}{{\Bbb Q}}

\def\1{\sqrt{-1}\:}

\newcommand{\cntrct}                
{\hspace{2pt}\raisebox{1pt}{\text{$\lrcorner$}}\hspace{2pt}}

\renewcommand{\phi}{\varphi}
\renewcommand{\epsilon}{\varepsilon}
\renewcommand{\geq}{\geqslant}
\renewcommand{\leq}{\leqslant}

\newcommand{\Teich}{\operatorname{\sf Teich}}

\newcommand{\Comp}{\operatorname{\sf Comp}}

\newcommand{\Kah}{\operatorname{Kah}}

\newcommand{\Pic}{\operatorname{Pic}}
\newcommand{\Pos}{\operatorname{Pos}}

\newcommand{\Aut}{\operatorname{Aut}}

\newcommand{\Mon}{\operatorname{\sf Mon}}

\newcommand{\Diff}{\operatorname{\sf Diff}}

\newcommand{\rk}{\operatorname{rk}}

\newcounter{Mycounter}[section]
\newcounter{lemma}[section]
\newcounter{claim}[section]
\renewcommand{\theclaim}{{Claim \thesection.\arabic{claim}}}
\newcommand{\claim}{%
    \setcounter{claim}{\value{Mycounter}}
    \refstepcounter{claim}
    \stepcounter{Mycounter}
    {\noindent \bf \theclaim:\ }}

\newcounter{sublemma}[section]
\newcounter{corollary}[section]
\renewcommand{\thecorollary}{{Corollary \thesection.\arabic{corollary}}}
\newcommand{\corollary}{%
    \setcounter{corollary}{\value{Mycounter}}
    \refstepcounter{corollary}
    \stepcounter{Mycounter}
    {\noindent \bf \thecorollary:\ }}

\newcounter{theorem}[section]
\renewcommand{\thetheorem}{{Theorem \thesection.\arabic{theorem}}}
\newcommand{\theorem}{%
    \setcounter{theorem}{\value{Mycounter}}
    \refstepcounter{theorem}
    \stepcounter{Mycounter}
    {\noindent \bf \thetheorem:\ }}

\newcounter{conjecture}[section]
\newcounter{proposition}[section]
\renewcommand{\theproposition}
      {{Proposition \thesection.\arabic{proposition}}}
\newcommand{\proposition}{%
    \setcounter{proposition}{\value{Mycounter}}
    \refstepcounter{proposition}
    \stepcounter{Mycounter}
    {\noindent \bf \theproposition:\ }}

\newcounter{definition}[section]
\renewcommand{\thedefinition}
      {{Definition~\thesection.\arabic{definition}}}
\newcommand{\definition}{%
    \setcounter{definition}{\value{Mycounter}}
    \refstepcounter{definition}
    \stepcounter{Mycounter}
    {\noindent \bf \thedefinition:\ }}

\newcounter{example}[section]
\newcounter{remark}[section]
\renewcommand{\theremark}{{Remark \thesection.\arabic{remark}}}
\newcommand{\remark}{%
    \setcounter{remark}{\value{Mycounter}}
    \refstepcounter{remark}
    \stepcounter{Mycounter}
    {\noindent \bf \theremark:\ }}

\newcounter{problem}[section]
\newcounter{question}[section]
\makeatletter

\@addtoreset{equation}{section} \@addtoreset{footnote}{section}
\makeatother

\def\blacksquare{\hbox{\vrule width 5pt height 5pt depth 0pt}}
\def\endproof{\blacksquare}

\begin{document}

\begin{center}
{\LARGE\bf
Construction of automorphisms of hyperk\"ahler manifolds
}

Ekaterina Amerik\footnote{Partially supported 
by RSCF grant 14-21-00053 within AG Laboratory NRU-HSE.}, Misha Verbitsky\footnote{Partially supported 
by RSCF grant 14-21-00053 within AG Laboratory NRU-HSE.


}

\end{center}

{\small \hspace{0.15\linewidth}
\begin{minipage}[t]{0.7\linewidth}
{\bf Abstract} \\
Let $M$ be an irreducible holomorphic symplectic (hyperk\"ahler) manifold. If $b_2(M)\geq 5$,
we construct a deformation $M'$ of $M$ which admits a symplectic automorphism 
of infinite order. This automorphism is hyperbolic, that is, its action on the space of real $(1,1)$-classes is hyperbolic. 
If $b_2(M) \geq 14$, similarly, we construct a deformation which admits a parabolic automorphism. 
\end{minipage}
}

\tableofcontents


\section{Introduction}


\subsection{Sublattices and automorphisms}

Since early 2000's, K3 surfaces were one of the prime subjects of
holomorphic dynamics (\cite{_Cantat:K3_}, \cite{_McMullen:K3_}).
By now, the dynamics of automorphism group acting on K3 surfaces is
pretty much understood (\cite{_Cantat_Dupont_}). Some of these
results are already generalized to irreducible holomorphic symplectic manifolds
of K\"ahler type (simple hyperk\"ahler manifolds) in any dimension
(\cite{_Oguiso:ICM_}).

The purpose of the present paper is to construct
sufficiently many interesting automorphisms on a
deformation of an arbitrary hyperk\"ahler manifold
(see Section \ref{_hype_intro_Section_} for basic definitions and properties of hyperk\"ahler manifolds).


For known examples of hyperk\"ahler manifolds, associated
with a K3 surface or an abelian surface, it is not
hard to find a deformation which admits a large automorphism group.
Indeed, we can lift an automorphism
of a K3 or a torus, or use some other explicit
construction. However, the classification problem for hyperk\"ahler manifolds still looks out of reach,
and finding deformations with interesting automorphism groups without referring to the explicit geometry
is much less obvious.

Even more complicated problem is to find $n$ to 1 rational correspondences (``rational
isogenies'') from a manifold to itself or to some
other hyperk\"ahler manifold. Such constructions are of considerable 
importance, but the visible ways to approach this
problem look rather difficult at the moment.

What makes possible the study of automorphisms, rather than isogenies, 
is that the group of automorphisms of a hyperk\"ahler
manifold can be understood in terms of its
period lattice  (that is, the Hodge structure
on the second cohomology and the BBF form, 
see Subsection \ref{BBF})
and the K\"ahler cone. The later is described
in terms of certain cohomology classes called
{\bf MBM classes} (\ref{_MBM_Definition_}, \ref{MBM-general}),
which are, roughly speaking, cohomology classes of negative BBF-square whose duals are represented by minimal rational
curves on a deformation of $M$.

This description is most easy to explain for a K3 surface.
In this case, MBM classes are  integral classes of
self-intersection $-2$, commonly called {\bf $(-2)$-classes}.

Let $M$ be a projective K3 surface, and $\Pos(M)\subset H^{1,1}(M)$
the positive cone, that is, the one of two connected
components of the set $\{v\in H^{1,1}(M,\R)\ \ |\ \ (v,v)>0\}$
which contains the K\"ahler classes. Denote by $R$ the set
of all $(-2)$-classes on $M$, and let $R^\bot$ be the union
of all orthogonal hyperplanes to all $v\in R$.
Then the K\"ahler cone $\Kah(M)$ is one of the 
connected components of $\Pos(M) \backslash R^\bot$, 
and 
\[
\Aut(M) = \{ g\in SO^+(H^2(M, \Z))\ \ |\ \ g(\Kah(M))=\Kah(M)\}.
\]
This gives an explicit description of the automorphism group, 
which becomes quite simple when $\Kah(M)=\Pos(M)$,
and this happens when $M$ has no $(-2)$-classes of Hodge type $(1,1)$.
When $\Kah(M)=\Pos(M)$, the group $\Aut(M)$ is identified
with the subgroup $\Gamma_M\subset SO^+(H^2(M, \Z))$,
\[
\Gamma_M=\{g\in SO^+(H^2(M,
\Z))\ \ |\ \ g(H^{1,1}(M))\subset  g(H^{1,1}(M))
\]
(``the group of Hodge isometries of $H^2(M,\Z)$''). It is not hard to see that $\Gamma_M$ is
mapped onto a finite index subgroup of the group of isometries of the Picard lattice
$\Pic(M)=H^{1,1}(M,\Z)$, hence it is infinite whenever this
lattice has infinite automorphism group. Now
$\Pic(M)$ has signature $(1,k)$ by Hodge index theorem. It is well-known
(see e.g. \cite{_Dickson_}) that $\Gamma_M$ is infinite
when $k>1$ and also when $k=1$ and the Picard lattice does not represent zero (that is, there is no nonzero 
$v\in \Pic(M)$ with $v^2=0$).

Therefore, to produce K3 surfaces with infinite
automorphism group, it would suffice to find a primitive
sublattice of rank $\geq 3$, signature $(1,k)$, $k\geq 2$
and without $(-2)$-vectors in $H^2(M,\Z)$. 
This can be done using the work of V. Nikulin, 
\cite{Nik}, which implies that any lattice
of signature $(1,k)$, $k<10$,  admits a primitive embedding to the K3 lattice (that is, an even
unimodular lattice of signature $(3, 19)$; such a lattice is unique up to isomorphism and isomorphic to $H^2(M,\Z)$).

The argument above produces {\bf symplectic  automorphisms},
that is, automorphisms which preserve the holomorphic
symplectic structure.

This approach is generalized in the present paper. 
In \cite{_AV:orbits_} it is shown that for each
hyperk\"ahler manifold $M$ there exists $N>0$, depending only on the deformation class of $M$, 
such that for all MBM classes $v$ one has $-N < q(v,v) <0$.
In the present paper, we prove that the lattice $H^2(M,\Z)$
of a hyperk\"ahler manifold $M$ satisfying $b_2(M)\geq 5$ (this is believed to hold always, 
but no proof exists today) 
contains a primitive sublattice $\Lambda\subset H^2(M,\Z)$ which does not represent numbers smaller than $N$
(that is, for any nonzero $v\in \Lambda$, one has $|q(v,v)|\geq N$).  .
Using the global Torelli theorem, we find 
a deformation $M_1$ of $M$ with $\Pic(M_1)=\Lambda$.
In this case, the Picard lattice of $M_1$ contains
no MBM classes, the K\"ahler cone coincides with the
positive cone, and the symplectic automorphism group is mapped
onto a finite index subgroup of the isometry group 
$O(\Lambda)$  (\ref{_Aut=Mon_Corollary_}).

This allows us to prove the following theorem.

\hfill

\theorem\label{_Automo_pic_rk_2_intro_Theorem_}
Let $M$ be a hyperk\"ahler manifold with 
$b_2(M)\geq 5$. Then $M$ admits a projective deformation
$M'$ with infinite group of symplectic automorphisms
and Picard rank 2.

{\bf Proof:} \ref{_hype_rank_2_Corollary_} \endproof

\hfill

The automorphisms obtained in
\ref{_Automo_pic_rk_2_intro_Theorem_}
are {\bf hyperbolic}: they act on $H^{1,1}(M)$
with one real eigenvalue $\alpha >1$,
another $\alpha^{-1}$, and the rest of eigenvalues
have absolute value 1. In fact, the symplectic automorphisms
of hyperk\"ahler manifolds can be classified in the same
way as automorphism of the hyperbolic plane
(see \ref{_classi_of_iso_hype_Theorem_}).
There are hyperbolic, or, more precisely, loxodromic 
automorphisms (ones which act on $H^{1,1}(M)$ with two real eigenvalues
of absolute value $\neq 1$), elliptic 
ones (automorphisms of finite order) and 
{\bf parabolic} (quasiunipotent with a non-trivial 
rank 3 Jordan cell).

If we want to produce a parabolic automorphism of a
deformation of a given hyperk\"ahler manifold, more
work is necessary. We need to find a primitive sublattice
$\Lambda \subset H^2(M, \Z)$ of signature $(1, k)$, $k\geq
2$, such that $q(v,v) \not\in ]-N, 0[$ for $v\in
    \Lambda$, and $\Lambda$ admits a parabolic isometry.
In order to produce such a sublattice we rely on the classification of 
rational vector spaces with a quadratic form by the signature, discriminant
and the collection of $p$-adic invariants, and on Nikulin's work on lattice
embeddings. Our method works under a stronger restriction on $b_2$.

The main problem (and the main reason for the 
strong restriction on $b_2$) is that the second cohomology lattice
$H^2(M,\Z)$ of a hyperk\"ahler manifold $M$ is not necessarily 
unimodular. In this case, one cannot apply Nikulin's
theorem directly. To construct the sublattice we need, 
we first embed our lattice $H^2(M,\Z)$ into
$H_\Q:=H\otimes_\Z \Q$, 
where $H$ is a unimodular lattice. Then we apply Nikulin's theorem to 
$H$, obtaining a primitive sublattice $\Lambda\subset H$, 
and take the intersection of $\Lambda$ with the image of
$H^2(M,\Z)\subset H_\Q$. This is no longer primitive in $H^2(M,\Z)$, but we have a good
control over the extent to which it is not, sufficient to assure that the ``primitivization''
does not have vectors of small nonzero square.

\hfill

\theorem
Let $M$ be a hyperk\"ahler 
manifold with $b_2(M)\geq 14$. Then $M$ has a deformation
with $\rk \Pic(M) \geq 3$, such that its group of symplectic
automorphisms contains a parabolic element.

{\bf Proof:} \ref{_pic_rk_14_parabo_Corollary_}.
\endproof

\hfill

\remark
Using the main result of \cite{_Verbitsky:ergodic_}, one sees that under the conditions of each of the two theorems, 
the points corresponding to hyperk\"ahler 
manifolds with a hyperbolic resp. parabolic automorphism are dense in the Teichm\"uller space.





\section{Hyperk\"ahler manifolds: basic results}
\label{_hype_intro_Section_}


In this section, we recall the definitions and basic
properties of hyperk\"ahler manifolds and MBM classes.

\subsection{Hyperk\"ahler manifolds}\label{BBF}

\definition
A {\bf hyperk\"ahler manifold} $M$, that is, a 
compact K\"ahler holomorphically symplectic manifold,
is called {\bf simple}, or {\bf maximal holonomy} hyperk\"ahler manifold 
(alternatively, {\bf irreducible holomorphically symplectic (IHS)}), if 
$\pi_1(M)=0$ and $H^{2,0}(M)=\C$.

\hfill

This definition is motivated by the following theorem
of Bogomolov. 

\hfill

\theorem \label{_Bogo_deco_Theorem_}
(\cite{_Bogomolov:decompo_})
Any hyperk\"ahler manifold admits a finite covering
which is a product of a torus and several 
simple hyperk\"ahler manifolds.
\endproof

\hfill

\remark
Further on, we shall tacitly assume that the hyperk\"ahler
manifolds we consider are of maximal holonomy (simple, IHS).

\hfill

The second cohomology $H^2(M,\Z)$ of a simple 
hyperk\"ahler manifold $M$ carries a
primitive integral quadratic form $q$, 
called {\bf the Bogomolov-Beauville-Fujiki form}. It generalizes the intersection product on a K3 surface:
its signature is $(3,b_2-3)$ on $H^2(M,\R)$ and $(1,b_2-3)$ on $H^{1,1}_{\R}(M)$.  It was first
defined in \cite{_Bogomolov:defo_} and 
\cite{_Beauville_},
but it is easiest to describe it using the
Fujiki theorem, proved in \cite{_Fujiki:HK_} and stressing the topological nature of the form.

\hfill

\theorem\label{_Fujiki_Theorem_}
(Fujiki)
Let $M$ be a simple hyperk\"ahler manifold,
$\eta\in H^2(M)$, and $n=\frac 1 2 \dim M$. 
Then $\int_M \eta^{2n}=c q(\eta,\eta)^n$,
where $q$ is a primitive integral nondegenerate quadratic form on $H^2(M,\Z)$, and 
$c>0$ is a rational number depending only on $M$. \endproof

\hfill






Consider $M$ as a differentiable manifold and denote by $I$ our complex structure on $M$ (we shall use notations like
$\Pic(M,I)$, $H^{1,1}(M,I)$ etc., to stress that we are working with this particular complex structure).
We call {\bf the Teichm\"uller space} $\Teich$ the quotient
$\Comp(M)/\Diff_0(M)$, where $\Comp(M)$
denotes the space of all complex structures of K\"ahler
type on $M$ and $\Diff_0(M)$ is the group of isotopies. It
follows from a result of Huybrechts (see
\cite{_Huybrechts:finiteness_}) that for an IHSM $M$,
$\Teich$ has only finitely many connected components. Let
$\Teich_I$ denote the one containing our given complex
structure $I$. Consider the subgroup of the mapping class
group $\Diff(M)/\Diff_0(M)$ fixing $\Teich_I$.

\hfill

\definition\label{monodr} The {\bf monodromy group}
$\Mon(M)$ is the image of this subgroup in $O(H^2(M,\Z), q)$. 
The {\bf Hodge monodromy group} $\Mon_I(M)$
is the subgroup of $\Mon(M)$ preserving the Hodge decomposition.

\hfill

\theorem\label{_Mon_arithmetic_Theorem_} 
(\cite{_V:Torelli_}, Theorem 3.5) 
The monodromy group is a finite index subgroup in $O(H^2(M, \Z), q)$.

\hfill

The image of the Hodge monodromy is therefore
an arithmetic subgroup of the orthogonal group of the 
Picard lattice $\Pic(M, I)$. Notice that the action of $\Mon_I(M)$ on the Picard lattice
can have a kernel; when $(M,I)$ is projective, it is easy to see that the kernel is a
finite group (just use the fact that it fixes a K\"ahler class and therefore consists of
isometries), but it can be infinite in general (\cite{_McMullen:K3_}).
By a slight abuse of notation, we sometimes also call the Hodge monodromy this
arithmetic subgroup itself; one way to avoid such an abuse is to introduce the 
{\bf symplectic Hodge monodromy group} $\Mon_{I, \Omega}(M)$ which is a subgroup
of $\Mon_I(M)$ fixing the symplectic form $\Omega$. Its representation on the Picard lattice
is faithful and the image is the same as that of $\Mon_I(M)$, so that $\Mon_{I, \Omega}(M)$
is identified to an arithmetic subgroup of $O(\Pic(M, I), q)$.

\hfill

\theorem (Markman's Hodge-theoretic Torelli theorem, \cite{_Markman:survey_})
The image of $\Aut(M, I)$ acting on $H^2(M)$ is the subgroup of $\Mon_I(M)$ preserving the
K\"ahler cone $\Kah(M, I)$.

\hfill



In this paper, we construct hyperk\"ahler manifolds
with the K\"ahler cone equal to the positive cone, 
and use this construction to find manifolds admitting interesting 
automorphisms of infinite order.

\subsection{MBM classes}

We call a cohomology class $\eta\in H^2(M, \R)$ 
{\bf positive} if $q(\eta,\eta)>0$, and 
{\bf negative} if $q(\eta,\eta)<0$. The {\bf positive cone} $\Pos(M,I)\in H^{1,1}_{\R}(M,I)$
is the connected component of the set 
of positive classes on $M$ which contains the K\"ahler classes. The K\"ahler cone is cut out inside the
positive cone by a certain, possibly infinite, number of rational hyperplanes (by a result of Huybrechts, 
we may take for these the orthogonals to the classes of rational curves).

\hfill

In \cite{_AV:MBM_}, we have introduced the following notion.

\hfill

\definition \label{_MBM_Definition_}
An integral $(1,1)$-class $z$ on $(M,I)$ is called {\bf monodromy birationally minimal (MBM)}, if for some
$\gamma\in \Mon_I(M)$, the hyperplane 
$\gamma(z)^{\bot}$ supports a (maximal-dimensional) face of the K\"ahler cone of a birational model of $(M,I)$.

\hfill

We have shown in \cite{_AV:MBM_} the invariance of the MBM property under all deformations of complex structure which leave $z$ of type $(1,1)$.
Moreover we have observed that a negative class $z$ generating the Picard group $\Pic(M,I)$ is MBM if and only if a 
rational 
multiple of $z$ is represented by a curve (in fact automatically rational; when we speak about curves representing
$(1,1)$-classes in cohomology, it means that we identify the integral classes of curves to certain rational 
$(1,1)$-classes
by the obvious isomorphism provided by the BBF form).
This leads to a simple extension of the notion to the
classes in the whole $H^2(M,\Z)$ rather than the Picard lattice.
By writing $M$ rather than $(M,I)$, we let a complex
structure $I$ vary in its deformation class; this class is
not uniquely determined by the topology, but there are
finitely many of them by the already-mentioned finiteness result of Huybrechts
(\cite{_Huybrechts:finiteness_}).

\



\definition\label{MBM-general}
A negative class $z\in H^2(M,\Z)$ on a hyperk\"ahler manifold
is called {\bf an MBM class} if there exist a deformation of $M$ 
with $\Pic(M)= \langle z \rangle$ such that $\lambda z$ is 
represented by a curve, for some $\lambda \neq 0$.

\hfill

\theorem\label{_MBM_faces_Theorem_}(\cite{_AV:MBM_}, Section 6)
Let $(M,I)$ be a hyperk\"ahler manifold, 
and $S$ the set of all its MBM classes of type $(1,1)$. The K\"ahler cone of $(M,I)$ is 
a connected component of $\Pos(M,I)\backslash \cup_{z\in S} z^\bot$.

\hfill

\remark 
As follows from an observation by Markman, the other connected components (``the K\"ahler chambers'') 
are the monodromy transforms of the K\"ahler cones
of birational models of $(M,I)$. The Hodge monodromy group permutes the K\"ahler chambers.

\hfill

From \ref{_MBM_faces_Theorem_} and Hodge-theoretic Torelli we easily deduce the following

\hfill

\corollary\label{_Aut=Mon_Corollary_}
Let $(M,I)$ be a hyperk\"ahler manifold which has no MBM classes
of type (1,1). Then any element of $\Mon_I(M)$ lifts to an automorphism of $(M,I)$.

\hfill

{\bf Proof:} Indeed, for such manifolds $\Kah(M,I)=\Pos(M,I)$ and therefore the whole group $\Mon_I(M)$ preserves 
the K\"ahler cone.
\endproof

\subsection{Morrison-Kawamata cone conjecture, MBM bound and automorphisms}

The following theorem has been proved in \cite{_AV:Mor_Kaw_}.

\hfill

\theorem (\cite{_AV:Mor_Kaw_}) Suppose that $(M,I)$ is projective and the Picard number $\rho(M, I)>3$. Then
the Hodge monodromy group has only finitely many orbits on the set of MBM classes
of type $(1,1)$ on $M$.
\endproof

\hfill

This result is a version of Morrison-Kawamata cone
conjecture for hyperk\"ahler manifolds. Its proof is based
on ideas of homogeneous dynamics (Ratner theory,
Dani-Margulis, Mozes-Shah theorems).

\hfill

Since the Hodge monodromy group acts 
by isometries, it follows that the primitive MBM classes in $H^{1,1}(M)$
have bounded square. Using deformations, one actually obtains the boundedness without the projectivity assumption
and with the condition $\rho(M, I)>3$ replaced by $b_2(M)>5$. One of the main tools of this paper
is a subsequent generalization of this statement.

\hfill

\theorem (\cite{_AV:orbits_})
Let $M$ be a hyperk\"ahler manifold with $b_2\geq 5$.
Then there exists a number $N>0$,
called  {\bf the MBM bound}, such that
any MBM class $z$ satisfies $|q(z,z)| <N$.

\hfill

Let us explain how this theorem permits 
one to construct hyperk\"ahler manifolds with large automorphism groups. 

\hfill

\definition
A {\bf lattice},  or a {\bf  quadratic lattice}, is a free abelian group 
$\Lambda\cong\Z^n$ equipped with an
integer-valued quadratic form $q$.  When we speak of an embedding
of lattices, we always assume that it is compatible
with their quadratic forms.

\hfill

\definition
A sublattice $\Lambda'\subset \Lambda$
is called {\bf primitive} if $\Lambda/\Lambda'$ is torsion-free.
A number $a$ is {\bf represented} by a lattice $(\Lambda,q)$ if $a= q(x,x)$
for some nonzero $x\in \Lambda$. 

\hfill

Now let $M$ be a hyperk\"ahler manifold. Consider the lattice $H^2(M,\Z)$ equipped with the BBF form $q$.
By Torelli theorem, for any primitive
sublattice $\Lambda\subset H^2(M,\Z)$ of signature $(1,k)$, there exists
a complex structure $I$ such that $\Lambda=\Pic(M, I)$ is the Picard lattice of $(M,I)$. The key remark is that as soon as we succeed in finding
such a primitive sublattice which does not represent small nonzero numbers, the corresponding hyperk\"ahler manifold 
has fairly large automorphism group. 

\hfill

\theorem\label{_aut_on_MBM_bound_Theorem_}
Let $M$ be a hyperk\"ahler manifold, and
$\Lambda\subset H^2(M,\Z)$ a primitive sublattice of signature $(1,k)$ which does not represent any number $a, 0\leq |a|\leq N$,
where $N$ is the MBM bound (we sometimes say in this case that $\Lambda$ ``satisfies the MBM bound''). Let $(M,I)$ be a deformation
of $M$ such that $\Lambda= \Pic(M, I)$. Then the K\"ahler cone of $(M,I)$ is equal to the positive cone and
the group of holomorphic symplectic automorphisms
$\Aut(M,\Omega)$ is projected with
finite kernel to $\Mon_{I, \Omega}(M)$, which is 
a subgroup of finite index in  $O(\Lambda)$.

\hfill

{\bf Proof:}
For the finiteness of the kernel of the natural map from $\Aut(M,\Omega)$ to 
$\Mon_{I, \Omega}(M)\subset GL(H^2(M))$ see e. g. \cite{_V:Torelli_}.
Since $\Lambda= H^{1,1}_I(M,\Z)$
satisfies the MBM bound, it contains no MBM classes, so the K\"ahler cone is equal to
the positive cone.
By \ref{_Aut=Mon_Corollary_}, $\Aut(M,\Omega)$ maps onto $\Mon_{I,\Omega}(M)$.
Now, $\Mon_{I,\Omega}(M)$ is a finite index subgroup in
$O(\Lambda)$, as follows from \ref{_Mon_arithmetic_Theorem_}.
\endproof


\section{Sublattices and automorphisms}


\subsection{Classification of automorphisms of a hyperbolic space}

\remark The group $O(m,n), m, n>0$ has 4 connected components.
We denote the connected component of 1 by $SO^+(m,n)$.
We call a vector $v$ {\bf  positive} if its square is positive.

\hfill

\definition
Let $V$ be a vector space with a quadratic form $q$ of signature 
$(1,n)$, $\Pos(V)=\{x\in V\ \ |\ \ q(x,x)>0\}$ 
its {\bf positive cone},  and ${\Bbb P}^+ V$ the projectivization 
of $\Pos(V)$.
Denote by $g$ any $SO(V)$-invariant Riemannian structure on
 ${\Bbb P}^+ V$. Then $({\Bbb P}^+ V, g)$ is called
{\bf hyperbolic space}, and the group $SO^+(V)$
{\bf the group of oriented  hyperbolic isometries}.

\hfill

\remark Since the isotropy group (the stabilizer of a point
$x\in {\Bbb P}^+ V$ for $SO^+(V)$-action on ${\Bbb P}^+ V$) 
is $SO(n)$ acting on $T_x {\Bbb P}^+ V= x^\bot$, the
hyperbolic metric on ${\Bbb P}^+ V$ is unique up to a
constant.

\hfill

\theorem \label{_classi_of_iso_hype_Theorem_}
{\bf (Classification of isometries of ${\Bbb P}^+ V$)}\\
Let $n>0$, and $\alpha \in SO^+(1,n)$ is an isometry
acting on $V$. Then one and only one of these three cases occurs
\begin{description}
\item[(i)]  $\alpha$ has an eigenvector $x$ with $q(x,x)>0$ 
($\alpha$ is ``an elliptic isometry'')
\item[(ii)] $\alpha$ has an eigenvector $x$ with $q(x,x)=0$ 
and real eigenvalue $\lambda_x$ satisfying $|\lambda_x|>1$
($\alpha$ is ``hyperbolic isometry'').
\item[(iii)] $\alpha$ has a unique eigenvector $x$ with
$q(x,x)=0$ and eigenvalue 1, and no fixed points on
${\Bbb P}^+ V$ ($\alpha$ is ``parabolic isometry'').
\end{description}

{\bf Proof:} This is a standard textbook result; see,
for instance, \cite{_Kapovich:Kleinian_}. \endproof

\hfill

\definition
Recall that the BBF form has signature $(1, b_2-3)$
on $H^{1,1}(M)$. An automorphism of a hyperk\"ahler manifold $(M,I)$
is called {\bf elliptic (parabolic, hyperbolic)}
if it is elliptic (parabolic, hyperbolic) on $H^{1,1}_I(M,\R)$.

\hfill




\subsection{Rank-two sublattices and existence of hyperbolic automorphisms}

In this Section, we prove the following theorem.

\hfill

\theorem\label{_sublattices_main_Theorem_hyperb} 
Let $L$ be a non-degenerate indefinite lattice of rank $\geq 5$, and $N$ a natural number.
Then $L$ contains a primitive rank 2 sublattice $\Lambda$ of signature $(1,1)$ which does not
represent numbers of absolute value less than $N$.

\hfill

This theorem immediately gives examples of hyperk\"ahler manifolds with hyperbolic automorphisms.

\hfill

\corollary\label{_hype_rank_2_Corollary_}
Let $M$ be a hyperk\"ahler 
manifold with $b_2(M)\geq 5$. Then $M$ has a deformation
admitting a hyperbolic automorphism.

\hfill

{\bf Proof:} Consider the lattice $L=H^2(M,\Z)$ and let $N$ be the
MBM bound for deformations of $M$. Take a sublattice $\Lambda$ as in \ref{_sublattices_main_Theorem_hyperb}  
and a deformation of $M$ such that $\Lambda= H^{1,1}_I(M,\Z)$.
Up to a finite index (meaning that the natural maps between these groups have finite kernel and image of finite index), 
$\Aut(M)=\Mon_I(M)=O(\Lambda)$. But $\Lambda$ does not represent zero, and then
it is well-known that $O(\Lambda)$ has a hyperbolic element (one way to view this is to interpret $\Lambda$, up to
a finite index, as 
a ring of integers in a real quadratic extension of $\Q$, and notice that the units provide automorphisms; so there
is an automorphism of infinite order, and it must automatically be hyperbolic as it cannot be parabolic or elliptic).
\endproof

\hfill

To prove \ref{_sublattices_main_Theorem_hyperb}, we need the following proposition.

\hfill

\proposition\label{_diag_nonzero-min_Proposition_}
Let $\Lambda$ be a diagonal rank 2 lattice with 
diagonal entries $\alpha_1, \alpha_2$ divisible by an 
odd power of $p$, $\alpha_i= \beta_i p^{2n_i+1}$, 
 and such that the numbers $\beta_i$ are not divisible by
$p$ and the equation $\beta_1 x^2 + \beta_2 y^2=0$ has no
solutions modulo $p$. Let $v\in \Lambda\otimes \Q$ be such that the value of the quadratic form 
on $v$ is an integer. Then this integer is divisible by $p$.

\hfill

{\bf Proof:} A direct computation, which is especially straightforward when one works in $\Q_p$
instead of $\Q$. \endproof

\hfill

{\bf Proof of \ref{_sublattices_main_Theorem_hyperb}}

\hfill

 By Meyer's Theorem \cite{_Meyers_},
$L$ has an isotropic vector
(that is, a vector $v$ with $q(v)=0$).
The {\bf isotropic quadric}
$\{v\in L \ \ |\ \ q(v)=0\}$ has infinitely
many points if it has one, and not all of them are proportional. Take two of such non-proportional points
$v$ and $v'$, and let $v_1:= av + b v'$.
Then $q(v_1)=2abq(v,v')$. We may chose $2ab$ to be of any
sign and such that it has arbitrary large 
prime divisors in odd powers. Concretely, consider the lattice $M=\langle v, v'\rangle^\bot$ of signature $(r-1,s-1)$
(here $(r,s)$ denotes the signature of $L$).
It is always possible to find a vector
$w\in \langle v, v'\rangle^\bot$ such that $q(w)$ is divisible
by an odd power of a suitable sufficiently large prime number $p$, but not by an even one (for instance consider a rank-two
sublattice where the form $q$ is equivalent to $x^2-dy^2$ over $\Q$ and pick a large $p$ such that $d$ is not a square
modulo $p$; then one can choose a suitable $w$ in such a sublattice). Now choose the multipliers $a, b$ in such a way that
the lattice $\Lambda:=\langle v_1, w\rangle$ satisfies 
assumptions of \ref{_diag_nonzero-min_Proposition_} with this $p$ and has signature $(1,1)$.
\endproof




\subsection{Sublattices of large rank and existence of parabolic automorphisms}

The purpose of this section is to construct, in $H^2(M,\Z)$, primitive sublattices of larger rank not 
representing small numbers
(except possibly zero). We use Nikulin's theorem on primitive embeddings into unimodular lattices.
As $H^2(M,\Z)$ is not always unimodular, we have to provide a trick which allows a reduction to the unimodular 
case. The trick consists in remarking that we can embed $H^2(M,\Z)$ into a lattice of rank
$b_2(M)+3$ which over $\Q$ is equivalent to a ``standard'' lattice $L_{st}=\sum \pm z_i^2$ (\ref{plus-3-dim}).
Then we take a suitable $\Lambda$ not representing numbers of absolute value less than $N$
(except possibly zero) and embed it into $L_{st}$ using Nikulin's results from \cite{Nik}.
The intersection $\Lambda\cap H^2(M,\Z)$ is not necessarily primitive in $H^2(M,\Z)$, but we can control
the extent to which it is not primitive in terms depending only on the embedding of $H^2(M,\Z)$ into $L_{st}$, 
and not on $N$; so, increasing $N$ if necessary, we eventually get a primitive sublattice satisfying 
the MBM bound.

\hfill

For a lattice $\Lambda$, we sometimes 
denote $\Lambda\otimes_{\Z}{\Q}$ by $\Lambda_{\Q}$.
Recall that the {\bf Hilbert symbol} $(a,b)_p$ of two 
$p$-adic numbers is equal to 1 if the equation
$ax^2+by^2=z^2$ has nonzero solutions in $\Q_p$ and -1 otherwise.
If $a$ and $b$ are nonzero rational numbers,
one has $(a,b)_p=1$ for all
$p$ except finitely many, and 
$\prod_p (a, b)_p=1$ (\cite{Se} chapter III, theorem 3).

\hfill

\theorem\label{plus-3-dim} For any non-degenerate lattice $(H, q)$ there is an embedding of rational vector spaces with a quadratic form
$(H\otimes_{\Z} \Q,q)\subset (L_{\Q}, q_{st})$, where $q_{st}=\sum \pm z_i^2$, the rank of $L_{\Q}$ is equal to $rk(H)+3$, and the signature
of $L_{\Q}$ can be taken arbitrary among the possible ones $(r+3, s)$; $(r+2, s+1)$; $(r+1, s+2)$; $(r, s+3)$, where $(r,s)$ is
the signature of $H$.

\hfill

{\bf Proof:} The form $q$ diagonalizes over the rationals; let $a_1, \dots, a_n$, $n=r+s$, be its diagonal entries and $d=a_1\dots a_n$.
It is well-known that a rational quadratic form is determined by its signature, its discriminant $d$ as an element of $\Q^*/(\Q^*)^2$ and
its collection of {\bf $p$-adic signatures}
$\epsilon_p(q)=\Pi_{i<j}(a_i, a_j)_p$ for all primes $p$,
where $(a_i, a_j)_p$ is the Hilbert symbol (\cite{Se}, Chapter IV.2, Theorem 9, Theorem 7).
 Let $t$ be the number of desired negative
diagonal entries for $(L_{\Q}, q_{st})$. We wish to add three dimensions to $H$ to make the quadratic form equivalent to $q_{st}$, that is, we are
looking for three extra diagonal entries $b_0$, $b_1$, $b_2=(-1)^tdb_0b_1$ such that for all primes $p$,
$$\epsilon_p(q_{st})=\epsilon_p(q)(d, b_0b_1b_2)_p(b_0, b_1b_2)_p(b_1,b_2)_p.$$

Using obvious identities like $(x, -x)_p=1$ and $(x, y^2z)_p=(x,z)_p$, we see that this amounts to asking 
that the quantities $(b_0, (-1)^{t-1}d)_p(b_1, (-1)^{t-1}db_0)_p$ have prescribed signs for all $p$, and finally that 
$((-1)^{t-1}db_0,(-1)^{t-1}db_1)_p$ have prescribed signs for every $p$. 
The prescribed signs must satisfy the usual identities for Hilbert symbols: only finitely many are equal to $-1$ and
the product is equal to one.

 In other words, to prove the theorem, it is sufficient 
to produce two rational numbers $x,y$ with arbitrarily defined signs such that for any $p$, $(x,y)_p$ are equal to prescribed $\delta_p$
(satisfying the usual identities). 
To do this, we use the following claim:

\hfill

\claim (\cite{Se}, Chapter III, Theorem 4) Let $x\in \Q$. To find an $y\in \Q$ with $(x,y)_p=\delta_p$ for all $p$, it suffices to find for each $p$, an $y_p\in \Q_p$ with
$(x, y_p)=\delta_p$.

\hfill

Now take $x$ which is not a square modulo all $p$ such that $\delta_p=-1$, then it is easy to find $y_p$ using the explicit formulae for the
Hilbert symbol (\cite{Se}, Chapter III, Theorem 1).
\endproof 

\hfill

Let now $H$ be $H^2(M,\Z)$ of signature $(3, b_2-3)$, $L_{\Q}$ as in \ref{plus-3-dim}, say, of signature $(3, b_2)$, and let $L$ be the standard 
integral lattice in $L_{\Q}$. This is an odd unimodular lattice. By the work of Nikulin, a lattice $\Lambda$ of signature $(1,s)$ and the
discriminant group $\Lambda^{\vee}/\Lambda$ without 2-torsion has a primitive embedding into $L$ as soon as  $s\leq b_2$ and $2(1+s) =2\rk \Lambda < \rk L=b_2+3$ (in fact Nikulin only gives an embedding result for even
lattices, but the crucial place in the argument is the construction of a lattice $\Lambda'$ with appropriate 
signature and the discriminant
form opposite to that of $\Lambda$; then the direct sum $\Lambda \oplus \Lambda'$ has a unimodular overlattice
corresponding to a maximal isotropic subgroup in the direct sum of the discriminant groups. This construction is done by Nikulin 
also in the odd case (\cite{Nik}, section 1.16), and yields a unimodular overlattice of $\Lambda \oplus \Lambda'$ exactly in the
same way).

In particular we can primitively embed into $L$ a lattice $\Lambda$ of signature $(1, [b_2/2])$, without 2-torsion in the discriminant
and not representing numbers of small absolute value other than zero (this is achieved for instance by multiplying a unimodular lattice by a large prime). 
We are actually looking for a primitive sublattice in $H$, whereas the intersection $\Lambda\cap H$, of signature 
$(1,[b_2/2]-3)$,  is not necessarily primitive in $H$. But its non-primitivity can be controlled in terms of the embedding of $H$ into $L_{\Q}$
and thus does not depend on $\Lambda$. In this way we obtain the following

\hfill

\theorem\label{emb-lattices-HK} 
Let $M$ be a hyperk\"ahler
manifold 
and $N$ a natural number. Then $H:=H^2(M,\Z)$ contains
a primitive sublattice of signature $(1, [b_2/2]-3)$ which does not represent non-zero numbers of absolute value smaller than $N$. 
In particular, there is a deformation of $M$ 
of Picard rank $[b_2/2]-2$ which does not have any MBM classes. 

\hfill

{\bf Proof:} Consider the embedding from $H$ into $L_{\Q}$ of dimension $b_2+3$ and signature $(3, b_2)$ as in \ref{plus-3-dim}.
Set $d=|H/H\cap L|$. 
Then for any primitive sublattice $\Lambda$ of $L$, the group 
$(H\cap \Lambda_\Q)/(H\cap L\cap \Lambda_\Q)=(H\cap \Lambda)_{\Q})/(H\cap \Lambda)$ embeds into 
$H/H\cap L$ and so has cardinality at most $d$. 
By Nikulin's results, if we take a lattice $\Lambda$ of signature $(1, [b_2/2])$,
without two-torsion in the discriminant and not representing any non-zero number
of absolute value less than $d^2N$, it admits a primitive
embedding to  $L$. Then the ``primitivization''
of $\Lambda\cap H$ (that is, $H\cap \Lambda_\Q \subset H$),
is a lattice of signature $(1, [b_2/2]-3)$ not
representing nonzero numbers of absolute value less than
$N$. Taking as $N$ the MBM bound for our manifold $M$, we obtain deformations
with relatively large Picard number and
no MBM classes.
\endproof

\hfill


\corollary\label{_pic_rk_14_parabo_Corollary_}
Let $M$ be a hyperk\"ahler 
manifold with $b_2(M)\geq 14$. Then $M$ has a deformation
admitting a parabolic automorphism.

\hfill

{\bf Proof. Step 1:}
Let $\Lambda= H^{1,1}_I(M,\Z)$ be a primitive lattice
of corank 2 and signature $(1,n)$ in $H^2(M,\Z)$ satisfying the MBM bound. 
Then $\Aut(M,\Omega)$ has finite index in  $O(\Lambda)$.
It suffices to show that the
Lie group $O(\Lambda\otimes_\Z \R)$ contains a rational
unipotent subgroup $U$. Then $U\cap \Aut(M,\Omega)$
is Zariski dense in $U$ by another application of Borel
and Harish-Chandra, and all its elements are parabolic. 

\hfill

{\bf Step 2:} Suppose that there exists a rational
vector $v$ with $q(v,v)=0$, and let $P\subset O(\Lambda\otimes_\Z \R)$
be the stabilizer of $v$. This subgroup is clearly
rational and parabolic; its unipotent radical is 
the group $U$ which we require.

\hfill

{\bf Step 3:} Such a rational vector exists for
any indefinite lattice of rank $\geq 5$ by Meyer's theorem
\cite{_Meyers_}, therefore as soon as $[b_2(M)/2]-2\geq 5$, $\Lambda$ has parabolic elements in its orthogonal group.
\endproof

\hfill

{\bf Acknowledgements:}
We are grateful to Noam Elkies, David Speyer and other participants of
the Mathoverflow discussion \cite{_Mathoverflow:sublattices_} for their
interest and advice, to Marat Rovinsky for his help with Hilbert symbols, to Slava Nikulin for patiently explaining certain points
of his work and to Dmitri Panov for his help on hyperbolic
geometry.

\hfill

{

{\small
\noindent {\sc Ekaterina Amerik\\
{\sc Laboratory of Algebraic Geometry,\\
National Research University HSE,\\
Department of Mathematics, 7 Vavilova Str. Moscow, Russia,}\\
\tt  Ekaterina.Amerik@gmail.com}, also: \\
{\sc Universit\'e Paris-11,\\
Laboratoire de Math\'ematiques,\\
Campus d'Orsay, B\^atiment 425, 91405 Orsay, France}

\hfill

\noindent {\sc Misha Verbitsky\\
{\sc Laboratory of Algebraic Geometry,\\
National Research University HSE,\\
Department of Mathematics, 7 Vavilova Str. Moscow, Russia,}\\
\tt  verbit@mccme.ru}, also: \\
{\sc Universit\'e Libre de Bruxelles, CP 218,\\
Bd du Triomphe, 1050 Brussels, Belgium}
 }
}

\end{document}